\documentclass{ptms-conf-2017}

\usepackage[utf8]{inputenc}

\usepackage{latexsym}
\newcommand\cA{{\cal A}}

\newcommand\cF{{\cal F}}

\newcommand\cB{{\cal B}}

\newcommand\cN{{\cal N}}

\newcommand\cR{{\cal R}}

\newcommand\e{\epsilon}
\newcommand\ve{\varepsilon}

\newcommand\Er{\mbox{Err}}

\def\bbr{{\mathbb R}}

\def\text#1{\hbox{#1}}

\def\E{{\bf E}}
\def\P{{\bf P}}

\def\C{{\bf C}}

\def\H{{\bf H}}

\def\Chi{{\bf 1}}

\def\d{\mathrm{d}}
\def\build #1_#2{\mathrel{\mathop{\kern 0pt #1}\limits_{#2}}}

\newcommand{\wh}{\widehat}
\newcommand{\wt}{\widetilde}
\newcommand{\zs}[1]{{\mathchoice{#1}{#1}{\lower.25ex\hbox{$\scriptstyle#1$}}
{\lower0.25ex\hbox{$\scriptscriptstyle#1$}}}}

\begin{document}
\mainmatter 

\title{Improved nonparametric estimation of the drift in diffusion processes}
\titlerunning{Improved nonparametric estimation of the drift in diffusion processes}
%

\author{${}^1$ Evgeny Pchelintsev, \, ${}^1$ Svyatoslav Perelevskiy \and ${}^1$ Irina Makarova}
\authorrunning{E. Pchelintsev, S. Perelevskiy, I. Makarova}

\institute{
	${}^1$Tomsk State University, 36 Lenina avenue, Tomsk, 634050, Russian Federation; \\
		\path|evgen-pch@yandex.ru|, \path|slavaperelevskiy@mail.ru|, \path|star_irish@bk.ru| \\
	\url{http://en.tsu.ru/}
}

\maketitle

\begin{abstract}
In this paper, we consider the robust adaptive non parametric estimation problem for the drift coefficient in diffusion processes. 
An adaptive model selection procedure, based on the improved weighted least square estimates, is proposed.  
Sharp oracle inequalities for the robust risk have been obtained.

\end{abstract}


\section{Introduction}\label{sec:In}
Let $(\Omega, \cF, (\cF)_\zs{t\ge 0},\P)$ be a  filtred probability space  on which the following stochastic differential equation  is defined:
 \begin{equation}\label{sec:In.1}
\d y_\zs{t}=S(y_\zs{t})\,\d t +\d w_\zs{t}\,,\quad 0\leq t\leq T\,,
 \end{equation}
where $(w_\zs{t})_{t\ge 0}$ is ascalar standard  Wiener process, the initial value $y_\zs{0}$ is some given constant, and
$S(\cdot)$ is an unknown function.The problem is to estimate the function $ S(x)$, $x in[a,b]$, from the observations $(y_t)_{0 \leq t \leq T}$.
The  calibration problem for the model  \eqref{sec:In.1}  is important in various applications. In particular, it appears, when constructing optimal strategies for investor behavior in diffusion financial markets. It is known that the optimal strategy depends on unknown market parameters, in particular, on  unknown drift coefficient $S$. Therefore, in practical financial calculations it is necessary to use statistical estimates for the function $S$ which  are reliable on some fixed time interval $[0,T]$ \cite{KaSh}.
Earlier, the problem of non-asymptotic estimation of the parameters of diffusion processes was studied in \cite{Kut2004}. Here it was shown that many difficulties of asymptotic estimation of parameters for one-dimensional diffusion processes can be overcome by using a sequantial approach. It turns out that the theoretical analysis of successive estimates is simpler than the analysis of classical procedures. In particular, it is possible to calculate non-asymptotic bounds for quadratic risk. Owing to the use of a sequential approach, the problems of non-asymptotic estimation of parameters were studied  in \cite{GaKo1} for multidimensional diffusion processes and recently in \cite{GaKo2} for multidimensional continuous and discrete semimartingales. In \cite{KoPe2} a truncated sequential method for estimating the parameters of diffusion processes was developed.
Now about nonparametric estimation. A consistent approach to nonparametric criteria for minimax estimation of the drift coefficient in (ergodic) diffusion processes was developed in \cite{GaPe4}. In this article, sequential pointwise kernel estimates are considered. For such estimates, non-asymptotic upper bounds of the root-mean-square risk are obtained, and these estimates give the optimal convergence rate as $ T \to \infty $. 

This paper deals with the  estimating the unknown function $S(x)$, $a\le x\le b$, in the sense of the mean square risk
 \begin{equation}\label{sec:In.2}
\cR(\wh{S}_\zs{T},S)\,=\,\E_\zs{S}\|\wh{S}_\zs{T}-S\|^2\,,
\quad
\|S\|^2=\int^b_\zs{a}\,S^2(x)\d x\,,
\end{equation}
where $\wh{S}_\zs{T}$ is the estimate of $S$ by observations $(y_\zs{t})_\zs{0\le t\le T}$, $a<b$ are some real numbers. Here $\E_\zs{S}$ is the expectation with respect to  the distribution $\P_\zs{S}$ of the random process $(y_\zs{t})_{0\le t\le T}$  given the drift function $S$. 

The goal of this paper is to construct an adaptive estimate $S^*$ of the drift coefficient $S$ in \eqref{sec:In.1} and to show that the quadratic risk of  this estimate is less then the one of the estimate proposed in \cite{GaPe4}, i.e. we construct the improved estimate in the mean square acuracy sence. For this we use the improved estimation approach proposed in \cite{Pch2013} and \cite{KPP2014} for parametric regression models and recently developted in \cite{PchPerg2017} for  a nonparametric estimation problem. Moreover in this paper we consider the estimation problem in adaptive setting, i.e. when the regulary of $S$ is unknown. For this we use a model selection method proposed in \cite{GaPe6}. Such approach provides adaptive solution for the nonparametric estimation through oracle inequalities which give the nonparametric upper bound for the quadratic risk of estimate.

The rest of the paper is organized as follows. In section 2 we reduce the initial problem to an estimation problem in a discrete time nonparametric regression model. In section 3 we construct the improved weigted least square estimates. In section 4 the sharp nonasymptotic oracle inequality for quadratic risk of model selection procedure is given.

\section{Passage to a discrete time regression model} \label{sec: Dt}
To obtain a reliable estimate of the function $S$, it is necessary to impose on it certain conditions that are analogous to the periodicity of the deterministic signal in the white noise model \cite{IbHa}. One of the conditions sufficient for this purpose is the assumption that the process $(y_\zs{t})_\zs{t\ge 0}$ in \eqref{sec:In.1} returns to any neighborhood of each points $ x \in [a, b] $.  As in \cite{GaPe4} to get the ergodicity of the process \eqref{sec:In.1} we define the following functional class:
    \begin{align}\nonumber
\Sigma_\zs{L,N}=\{&S\in Lip_\zs{L}(\bbr)\,:\,|S(N)|\le L \,; \
 \forall |x|\ge N,\ \exists\ \dot{S}(x)\in \C(\bbr)
\\[2mm] \label{sec:Dt.1} 
& \mbox{\rm such that}
-L\le \inf_\zs{|x|\ge N} \dot{S}(x)\le
\sup_\zs{|x|\ge N}\dot{S}(x)\le-1/L\}\,,
 \end{align}
where  $L>1$, $N>|a|+|b|$, $\dot{S}(x) - $ derivative $S(x)$,
$$
Lip_\zs{L}(\bbr)\,=\,
\left\{f\in\C(\bbr)\,:\,\sup_\zs{x,y\in\bbr}\frac{|f(x)-f(y)|}{|x-y|}\,\le\,L\right\}\,.
$$
     
  We note that if $S\in \Sigma_\zs{L,N}$,
then there exists an invariant density
 \begin{equation}\label{sec:Dt.2}
q(x)\,=q_\zs{S}(x)\,=\,\frac{\exp\{2\int^x_0 S(z)\d z\}}{
\int^{+\infty}_{-\infty}\exp\{2\int^y_0 S(z)\d z\}\d y}.
 \end{equation}
 
We note that the functions in
$\Sigma_\zs{L,N}$  are uniformly bounded on  $[a,b]$, i.e.
\begin{equation*}\label{sec:Dt.5}
s^{*}=\sup_\zs{a\le x\le b}
\sup_\zs{S\in\Sigma_\zs{L,N}}
S^2(x)\,<\infty\,.
\end{equation*}

 We start with the partition of the interva $[a,b]$ by the points $(x_\zs{k})_\zs{1\le k\le n}$, defined as
\begin{equation}\label{sec:Dt.6}
x_\zs{k}=a+\frac{k}{n} (b-a)\,,
\end{equation}
where $n=n(T)$ is an integer-valued function of $T$
such that
\begin{equation}\label{sec:Dt.7}
n(T)\le T
\quad\mbox{and}\quad
\lim_\zs{T\to\infty}\frac{n(T)}{T}=1\,.
\end{equation}

\noindent Now at any point $x_\zs{k}$ we estimate the function  $S$ by a sequential kernel estimation. We fix some $0<t_\zs{0}<T$ and put

\begin{equation}\label{sec:Dt.8}
\left\{
\begin{array}{cl}
\tau_\zs{k}&=\inf\{t\ge t_\zs{0}\,:\,\int^t_\zs{t_\zs{0}}\,
Q\left(\frac{y_\zs{s}-x_\zs{k}}{h}\right)\,\d s
\ge\,H_\zs{k}\}\,;\\[6mm]
\wt{S}_\zs{k}&=\dfrac{1}{H_\zs{k}}\,
\int^{\tau_\zs{k}}_\zs{t_\zs{0}}\,Q\left(\frac{y_\zs{s}-x_\zs{k}}{h}\right)\,\d y_\zs{s}\,,
\end{array}
\right.
\end{equation}
where $Q(z)=\Chi_\zs{\{|z|\le 1\}}$, $\Chi_\zs{A}$ is an indicator of the set $A$,  $h=(b-a)/(2n)$ and $H_\zs{k}$ is a positive threshold, which will be indicated below.
From \eqref{sec:In.1}
it is easy to obtain that
$$
\wt{S}_\zs{k}=S(x_\zs{k})+\zeta_\zs{k}\,.
$$
The error $\zeta_\zs{k}$ is represented as a sum of the approximating and stochastic parts, i.e.
$$
\zeta_\zs{k}=B_\zs{k}+\frac{1}{\sqrt{H_\zs{k}}}\,\xi_\zs{k}\,,\quad
B_\zs{k}=\frac{1}{H_\zs{k}}\,
\int^{\tau_\zs{k}}_\zs{t_\zs{0}}\,Q\left(\frac{y_\zs{s}-x_\zs{k}}{h}\right)\,
\Delta S(y_\zs{s},x_\zs{k})\d s\,,
$$
where $\Delta S(y,x)=S(y)-S(x)$ and
$$
\xi_\zs{k}=\frac{1}{\sqrt{H_\zs{k}}}\,
\int^{\tau_\zs{k}}_\zs{t_\zs{0}}\,Q\left(\frac{y_\zs{s}-x_\zs{k}}{h}\right)\,\d w_\zs{s}\,.
$$
Taking into account that $S$   is Lipshitz function , we obtain an upper bound for the approximating part as
\begin{equation*}\label{sec:Dt.9}
|B_\zs{k}|\le L h\,.
\end{equation*}
It is easy to see that random variables
$(\xi_\zs{k})_\zs{1\le k\le n}$ are independent identically distributed from  $\cN(0,1)$.
In \cite{GaPe4} it is established that an effective kernel estimate of the form \eqref{sec:Dt.8} has a stochastic part distributed as
$\cN(0,2Th q_\zs{S}(x_\zs{k}))$, where $q_\zs{S}(x_\zs{k})$  is the ergodic density defined in \eqref{sec:Dt.2}.
Therefore, for an effective estimate at each point $x_\zs{k}$ by the kernel estimate \eqref{sec:Dt.8},
 we need to estimate the density \eqref{sec:Dt.2}
from observations $(y_\zs{t})_\zs{0\le t\le t_\zs{0}}$.To this end, we establish that
$$
\wt{q}_T(x_\zs{k})=\max\{\wh{q}(x_\zs{k})\,,\,\epsilon_\zs{T}\}\,,
$$
where $\epsilon_\zs{T}$ 
is positive, $0<\epsilon_\zs{T}<1$,
\begin{equation*}\label{sec:Dt.10}
\wh{q}(x_\zs{k})=
\frac{1}{2t_\zs{0}h}
\int^{t_\zs{0}}_\zs{0}Q\left(\frac{y_\zs{s}-x_\zs{k}}{h}\right)\d s\,.
\end{equation*}
Now choose the threshold $H_\zs{k}$ in \eqref{sec:Dt.8}:
$$
H_\zs{k}=(T-t_\zs{0})(2\wt{q}_T(x_\zs{k})-\epsilon^2_\zs{T})h\,.
$$
Suppose that the parameters
$t_\zs{0}=t_\zs{0}(T)$ and $\epsilon_\zs{T}$ satisfy the following conditions:
\begin{itemize}
\item[$\H_\zs{1})$]
{\sl For any $T\ge 32$,
$$
16\le t_\zs{0}\le T/2
\quad\mbox{and}\quad
\sqrt{2}/t^{1/8}_0\le\epsilon_\zs{T}\le 1\,.
$$
}
\item[$\H_\zs{2})$]
$$
\lim_\zs{T\to\infty}\,t_\zs{0}(T)\,=\,\infty\,,\quad
\lim_\zs{T\to\infty}\,\e_\zs{T}\,=\,0\,,\quad
\lim_\zs{T\to\infty}\,T\e_\zs{T}/t_\zs{0}(T)\,=\,\infty\,.
$$

\item[$\H_\zs{3})$]
{\sl For any $\nu>0$ and $m>0$,
$$
\lim_\zs{T\to\infty}T\e^m_\zs{T}=\infty
\quad\mbox{and}\quad
\lim_\zs{T\to\infty}\,
T^{m}\,e^{-\nu \sqrt{t_\zs{0}}}\,=0\,.
$$
}
\end{itemize}
For example, for $T\ge 32$,
$$
t_0=\max\{\min\{\ln^{4} T\,,\,T/2\}\,,\,16\}
\quad\mbox{and}\quad
\epsilon_\zs{T}=\sqrt{2}\,t^{-1/8}_0\,.
$$
Let
\begin{equation}\label{sec:Dt.11}
\Gamma=\{\max_\zs{1\le l\le n}\tau_\zs{l}\,\le\,T\}
\quad\mbox{and}\quad
Y_\zs{k}=\wt{S}_\zs{k}\,\Chi_\zs{\Gamma}\,.
\end{equation}
Then on the set $\Gamma$ there exists a temporary heteroscedastic regression model
\begin{equation}\label{sec:Dt.12}
Y_\zs{k}=S(x_\zs{k})+\zeta_\zs{k}\,,\quad \zeta_\zs{k}=\sigma_\zs{k}\,\xi_\zs{k}+\delta_\zs{k}
\end{equation}
with $\delta_\zs{k}=B_\zs{k}$ and
\begin{equation*}\label{sec:Dt.12-1}
\sigma^{2}_\zs{k}=\frac{n}{(T-t_\zs{0})(\wt{q}_T(x_\zs{k})-\epsilon^2_\zs{T}/2)(b-a)}\,.
\end{equation*}
\noindent It should be noted that from \eqref{sec:Dt.7} and $\H_\zs{1})$, we get the following upper bound
\begin{equation}\label{sec:Dt.14}
\max_\zs{1\le k\le n}\sigma^2_\zs{k}\le \frac{4}{(b-a)\epsilon_\zs{T}}
=\sigma_\zs{*}
\end{equation}
for which, by condition $\H_\zs{3})$,
\begin{equation*}\label{sec:Dt.14-1}
\lim_\zs{T\to\infty}
\frac{\sigma_\zs{*}}{T^m}=0
\quad\mbox{for any}\quad m>0\,.
\end{equation*}

\noindent To estimate the $ S $ function from the observations of  \eqref{sec:Dt.12}
   should study some properties of the set $\Gamma$ in \eqref{sec:Dt.11}.
\begin{proposition}\label{Pr.sec:Dt.1}
Suppose that the parameters $t_\zs{0}$ and $\epsilon_\zs{T}$ satisfy the following conditions:
$\H_\zs{1})$ -- $\H_\zs{3})$. Then
\begin{equation*}\label{sec:Dt.15}
\sup_\zs{S\in\Sigma_\zs{L,N}}
\P_\zs{S}(\Gamma^c)
\le\,
\Pi_\zs{T}\,,
\end{equation*}
where $\lim_\zs{T\to\infty}\,T^{m}\,\Pi_\zs{T}=0$ for any $m> 0$.
\end{proposition}

\section{Improved estimates}\label{sec:Ms}

In this section we consider the estimation problem for the model \eqref{sec:Dt.12}.
The function $S(\cdot)$ is unknown and has to be estimated
from observations $Y_1,\ldots,Y_n$.

The accuracy of any estimator $\wh{S}$ will be measured by
the empirical squared error of the form
$$
\|\wh{S}-S\|^2_n=(\wh{S}-S,\wh{S}-S)_\zs{n}
=\frac{b-a}{n}\sum^n_\zs{l=1}(\wh{S}(x_\zs{l})-S(x_\zs{l}))^2\,.
$$
Now we fix a basis 
$(\phi_\zs{j})_\zs{1\le j\le n}$ which is orthonormal for the empirical inner product:
\begin{equation*}\label{sec:Ms.4}
(\phi_i\,,\,\phi_\zs{j})_\zs{n}=
\frac{b-a}{n} \sum^n_\zs{l=1}\,\phi_i(x_\zs{l})\phi_\zs{j}(x_\zs{l})= {\bf Kr}_\zs{ij}\,,
\end{equation*}
where ${\bf Kr}_\zs{ij}$ is Kronecker's symbol.
By making use  of this basis we apply the discrete Fourier transformation to
\eqref{sec:Dt.12} and we obtain the Fourier coefficients
$$
\wh{\theta}_\zs{j,n}=\frac{b-a}{n}\sum^n_\zs{l=1}\,Y_\zs{l}\phi_\zs{j}(x_\zs{l})\,,
\quad
\theta_\zs{j,n}=\frac{b-a}{n}\sum^n_\zs{l=1}S(x_\zs{l})\,\phi_\zs{j}(x_\zs{l})\,.
$$
From \eqref{sec:Dt.12} it follows directly that these Fourier coefficients satisfy the following
equation
\begin{equation*}\label{sec:Ms.5}
\wh{\theta}_\zs{j,n}\,=\,\theta_\zs{j,n}\,+\,\zeta_\zs{j,n}
\quad\mbox{with}\quad
\zeta_\zs{j,n}=\sqrt{\frac{b-a}{n}}\xi_\zs{j,n}+\delta_\zs{j,n}\,,
\end{equation*}
where
$$
\xi_\zs{j,n}=\sqrt{\frac{b-a}{n}}\sum^n_\zs{l=1}\sigma_\zs{l}\xi_\zs{l}\phi_\zs{j}(x_\zs{l})
\quad\mbox{and}\quad
\delta_\zs{j,n}=\frac{b-a}{n}\sum^n_\zs{l=1}\,\delta_\zs{l}\,\phi_\zs{j}(x_\zs{l})\,.
$$
Note that the upper bound \eqref{sec:Dt.14} and the Bounyakovskii-Cauchy-Schwarz inequality imply that
\begin{equation*}\label{sec:Ms.6}
|\delta_\zs{j,n}|\le\|\delta\|_\zs{n}\,\|\phi_\zs{j}\|_\zs{n}=
\|\delta\|_\zs{n}\,.
\end{equation*}
We estimate the function $S$ in \eqref{sec:Dt.12} on the sieve 
\eqref{sec:Dt.6} by the weighted least squares estimator
\begin{equation*}\label{sec:Ms.7}
\wh{S}_\zs{\lambda}(x_\zs{l})\,=\,\sum^n_\zs{j=1}\,\lambda(j)\,\wh{\theta}_\zs{j,n}\,\phi_\zs{j}(x_\zs{l})\,
\Chi_\zs{\Gamma}\,,\quad 1\le l\le n\,,
\end{equation*}
where the weight vector $\lambda=(\lambda(1),\ldots,\lambda(n))$ 
belongs to some finite set $\Lambda\subset [0,1]^n$. 
We set for any $a\le x\le b$
\begin{equation}\label{sec:Ms.8}
\wh{S}_\zs{\lambda}(x)=
\wh{S}_\zs{\lambda}(x_\zs{1})
\Chi_\zs{\{a\le x\le x_\zs{1}\}}+
\sum^{n}_\zs{l=2}\wh{S}_\zs{\lambda}(x_\zs{l})
\Chi_\zs{\{x_\zs{l-1}< x\le x_\zs{l}\}}\,.
\end{equation}

Further we suppose that the first $d \le n$ components of the weight vector $\lambda$ are equal to 1, i.e. $\lambda(j)=1$ for any $1\le j\le d$.

We consider a new estimate for the function $S$ in \eqref{sec:Dt.12} of the form
\begin{equation*}\label{sec:Imp.1}
S_\zs{\lambda}^{*}(x_\zs{l})\,=\,\sum^n_\zs{j=1}\,\lambda(j)\,\theta_\zs{j,n}^{*}\,\phi_\zs{j}(x_\zs{l})\,
\Chi_\zs{\Gamma}\,,\quad 1\le l\le n\,,
\end{equation*}
where 
\begin{equation*}\label{sec:Imp.2}
\theta_\zs{j,n}^{*}=\left(1-\frac{c(d)}{\|\wt{\theta}_\zs{n}\|}\Chi_\zs{\{1\le j\le d\}}\right)\wh{\theta}_\zs{j,n},
\end{equation*}
 where  
 \begin{equation*}
 c(d)=\frac{(d-1)\sigma^2_*L(b-a)^{1/2}}{n(s^*+\sqrt{d\sigma_*/n})}\,,\quad \|\wt{\theta}_\zs{n}\|^2=\sum_{j=1}^{d}\wh{\theta}^2_\zs{j,n}.
 \end{equation*}
Now we define the estimate for $S$ in \eqref{sec:In.1}. We set for any $a\le x\le b$ 
\begin{equation}\label{sec:Imp.3}
S_\zs{\lambda}^{*}(x)=
S_\zs{\lambda}^{*}(x_\zs{1})
\Chi_\zs{\{a\le x\le x_\zs{1}\}}+
\sum^{n}_\zs{l=2}S_\zs{\lambda}^{*}(x_\zs{l})
\Chi_\zs{\{x_\zs{l-1}< x\le x_\zs{l}\}}\,.
\end{equation}

We denote the difference of quadratic risks of the estimates \eqref{sec:Imp.3} and \eqref{sec:Ms.8} as
$$
\Delta_n(S):=\E_\zs{S}\|S^{*}_\zs{\lambda}-S\|^2_\zs{n}-\E_\zs{S}\|\wh{S}_\zs{\lambda}-S\|^2_\zs{n}.
$$

The choice of estimate \eqref{sec:Imp.3} is motivated by the
desire to control the quadratic risk.

\begin{theorem}
The estimate \eqref{sec:Imp.3} outperforms in mean square accuracy the estimate \eqref{sec:Ms.8}, i.e.
$$
\sup_{S \in \Sigma_{L,N}}\Delta_n(S)<-c^2(d).
$$
\end{theorem}

\section{Oracle inequalities}\label{sec:Or}

 In order to obtain
a good estimator, we have to write a rule to choose a weight vector 
$\lambda\in\Lambda$ in \eqref{sec:Imp.3}. 
It is obvious, that the best way is to minimize 
 the empirical squared error with respect to $\lambda$:
$$
\Er_\zs{n}(\lambda)=\|S^*_\zs{\lambda}-S\|^2_\zs{n} \to\min\,.
$$
Making use of \eqref{sec:Imp.3} and the Fourier transformation of $S$ imply
\begin{equation*}\label{sec:Ms.10}
\Er_\zs{n}(\lambda)=
\sum^n_\zs{j=1}\,\lambda^2(j) \theta^{*2}_\zs{j,n}\,-
2\,\sum^n_\zs{j=1}\,\lambda(j) \theta^{*}_\zs{j,n}\,\theta_\zs{j,n}\,+\,
\sum^n_\zs{j=1}\,\theta^2_\zs{j,n}\,.
\end{equation*}
Since the coefficient $\theta_\zs{j,n}$ is unknown, we need to replace the term 
$\theta^{*}_\zs{j,n}\,\theta_\zs{j,n}$
by some its estimator which we choose as 
\begin{equation*}\label{sec:Ms.11}
\wt{\theta}_\zs{j,n}=
\wh{\theta}_\zs{j,n}\theta^{*}_\zs{j,n}-\frac{b-a}{n}s_\zs{j,n}
\quad\mbox{with}\quad
s_\zs{j,n}=\frac{b-a}{n}\,\sum^n_\zs{l=1}\,\sigma^2_\zs{l}\,\phi^2_\zs{j}(x_\zs{l})\,.
\end{equation*}
One has to pay a penalty for this substitution in the empirical squared error. 
Finally, we define the cost function of the form
\begin{equation*}\label{sec:Ms.12}
J_\zs{n}(\lambda)\,=\,\sum^n_\zs{j=1}\,\lambda^2(j)\theta^{*2}_\zs{j,n}\,-
2\,\sum^n_\zs{j=1}\,\lambda(j)\,\wt{\theta}_\zs{j,n}\,
+\,\rho P_\zs{n}(\lambda)\,,
\end{equation*}
where  the penalty term is defined as
$$
P_\zs{n}(\lambda)=
\frac{b-a}{n}\sum^n_\zs{j=1}\lambda^2(j)s_\zs{j,n}
$$
and $0< \rho< 1$ is some positive constant which will be chosen later.
We set
\begin{equation*}\label{sec:Ms.13}
\wh{\lambda}=\mbox{argmin}_\zs{\lambda\in\Lambda}\,J_n(\lambda)
\end{equation*}
and  define an estimator of $S$ of the form \eqref{sec:Ms.8}:
\begin{equation}\label{sec:Ms.14}
S^{*}(x)=S^*_\zs{\wh{\lambda}}(x)
\quad\mbox{for}\quad a\le x\le b\,.
\end{equation}

\noindent Now we obtain the non asymptotic upper bound for the quadratical risk 
of the estimator \eqref{sec:Ms.14}.
\begin{theorem}\label{Th.sec:Ms.1}
Let $\Lambda\subset [0,1]^n$ be any finite set such that  the first $d \le n$ components of the weight vector $\lambda$ are equal to 1. Then, for any $n\ge 3$ and $0<\rho<1/6$,
the estimator \eqref{sec:Ms.14} satisfies
 the following oracle inequality
$$
\E_\zs{S}\|S^{*}-S\|^2_\zs{n}\le
\frac{1+6\rho}{1-6\rho}\,
\min_\zs{\lambda\in\Lambda}\,
\E_\zs{S}\|\wh{S}_\zs{\lambda}-S\|^2_\zs{n}
+
\frac{\Psi_\zs{n}(\rho)}{n}\,,
$$
where
$\lim_{n\to\infty}\Psi_\zs{n}(\rho)/n=0$.
\end{theorem}

Now we consider the estimation problem \eqref{sec:In.1} via model 
\eqref{sec:Dt.12}.
We apply the estimating procedure 
\eqref{sec:Ms.14} with special weight  
set introduced in \cite{GaPe4} to the regression scheme \eqref{sec:Dt.12}.  Denoting $S^*_\zs{\alpha}=S^*_\zs{\lambda_\zs{\alpha}}$ we set
\begin{equation*}\label{sec:Ms.24}
S^{*}=S^*_\zs{\wh{\alpha}} \quad\mbox{with}\quad
\wh{\alpha}=\mbox{argmin}_\zs{\alpha\in\cA_\zs{\ve}}\,J_n(\lambda_\zs{\alpha})\,.
\end{equation*}

We obtain through
 Theorem~\ref{Th.sec:Ms.1} the following oracle inequality.
\begin{theorem}\label{Th.sec:Ms.2}
Assume that $S\in\Sigma_\zs{L,N}$ and the number of the points 
$n=n(T)$ in the model\eqref{sec:Dt.12} satisfies 
\eqref{sec:Dt.7}.
Then the procedure  $S^{*}$ satisfies, 
for any $T\ge 32$, the following inequality 
\begin{equation*}\label{sec:Ms.25}
\cR(S^{*},S)
\,\le\,\frac{(1+\rho)^2(1+6\rho)}{1-6\rho}
\,\min_\zs{\alpha\in\cA_\zs{\ve}}\,\cR(S^*_\zs{\alpha},S)
+
\frac{\cB_\zs{T}(\rho)}{n}\,,
\end{equation*}
where 
$
\lim_{T\to\infty}\cB_\zs{T}(\rho)/n(T)=0.
$
\end{theorem}

\section*{Acknowledgement}

The results   of Section 3 of this work are supported by the RSF grant number 17-11-01049. The results of section 4 are supported by the Ministry of Education and Science of the Russian Federation in the framework of the research project no. 2.3208.2017/4.6.


\end{document}